\documentclass[a4paper,12pt]{article}
\usepackage{a4wide}
\usepackage{theorem}
\usepackage{amsmath}
\usepackage{array}
\usepackage{amssymb}
\usepackage{amsfonts}
\usepackage[english]{babel}
\usepackage{graphicx}

\newtheorem{theo}{\indent Theorem \newline}[section]
\newtheorem{defi}[theo]{\indent Definition\newline}
\newtheorem{rem}[theo]{\noindent Remark}
\newtheorem{prop}[theo]{\indent Proposition\newline}
\newtheorem{lemma}[theo]{\indent Lemma\newline}
\newtheorem{cor}[theo]{\indent Corollary \newline}

 \def\N{{\mathbb{N}}}
\def\Z{{\mathbb{Z}}}

\def\R{{\mathbb{R}}}
\def\C{{\mathbb{C}}}
\def\H{{\mathbb{H}}}

\usepackage{eufrak}

\newlength{\indentation}%
\makeatletter
\newcommand\@makefntextsans[1]{%
    \parindent 0em%
    \noindent%
    \hb@xt@0em{\hss}%
    #1}
\def\footnotetextsans{%
     \@ifnextchar [\@xfootnotenextsans%
       {\@footnotetextsans}}
\def\@xfootnotenextsans[#1]{%
  \begingroup%
     \csname c@\@mpfn\endcsname #1\relax%
  \endgroup%
  \@footnotetextsans}
\long\def\@footnotetextsans#1{\insert\footins{%
    \reset@font\footnotesize%
    \interlinepenalty\interfootnotelinepenalty%
    \splittopskip\footnotesep%
    \splitmaxdepth \dp\strutbox \floatingpenalty \@MM%
    \hsize\columnwidth \@parboxrestore%
    \color@begingroup%
      \@makefntextsans{%
        \rule\z@\footnotesep\ignorespaces#1\@finalstrut\strutbox}
    \color@endgroup}}
\makeatother

\begin{document}

\cleardoublepage
\title{Open Gromov-Witten invariants in dimension six}
\author{Jean-Yves Welschinger}
\maketitle

\makeatletter\renewcommand{\@makefnmark}{}\makeatother
\footnotetextsans{Keywords: holomorphic discs, Gromov-Witten invariants.}
\footnotetextsans{AMS Classification : 53D45.
}

{\bf Abstract:}

Let $L$ be a closed orientable Lagrangian submanifold of a closed symplectic six-manifold $(X , \omega)$.
We assume that the first homology group $H_1 (L ; A)$ with coefficients in a commutative ring $A$ injects into the group
$H_1 (X ; A)$ and that $X$ contains no Maslov zero pseudo-holomorphic disc with boundary on $L$. 
Then, we prove that for every generic choice of a tame almost-complex structure $J$ on $X$, every relative homology class $d \in H_2 (X , L ; \Z)$
and adequate number of incidence conditions in $L$ or $X$,
the weighted number of $J$-holomorphic discs with boundary on $L$, homologous to $d$, and either irreducible or reducible disconnected,
which satisfy the conditions,  does not depend on the generic choice of $J$, provided that at least one incidence condition lies in $L$. 
These numbers thus define open Gromov-Witten invariants in dimension six, taking values in the ring $A$. 

\section*{Introduction}
 
 In the mid eighties, M. Gromov discovered that classical enumerative invariants of complex  geometry
 obtained by counting the number of curves satisfying some incidence conditions in a smooth projective manifold
 actually only depend on the underlying K\"ahler form of the manifold up to deformation and not that much on the
 algebraic structure. A famous example of such an enumerative invariant is the number of  degree $d$ rational curves 
 passing through $3d-1$ generic points in the complex projective plane, a number later computed by M. Kontsevich. 
 The strategy followed by Gromov was to first introduce an auxiliary generic almost-complex structure tamed by the 
 K\"ahler form and then to count in an appropriate way the finite number of $J$-holomorphic curves satisfying the
 incidence conditions. He could then prove that this number does not depend on the generic choice of the almost-complex
 structure $J$. This approach and these results gave birth to the still developing theory of Gromov-Witten invariants in
 symplectic geometry.  A question then appeared together with these works of Gromov and (later) Witten, \cite{Gro}, \cite{Wit}. 
 Given a closed Lagrangian submanifold $L$ of a closed symplectic manifold $(X , \omega)$, is it likewise possible to extract
 enumerative invariants from the count of $J$-holomorphic discs of $X$ with boundary on $L$ and subject to given incidence conditions,
 in the sense that this number does not depend on the generic choice of $J$? Though apparently similar, this question hides
 a new difficulty, namely, the moduli spaces of $J$-holomorphic discs have real codimension one boundary components,
 contrary to the moduli spaces of closed $J$-holomorphic curves. As a consequence, even the invariance modulo two does not hold in general. 
 
  At the early 2000's, such open Gromov-Witten invariants have been defined by C.C. Liu and M. Katz in the presence of an action of the circle, 
  see \cite{KatzLiu}, \cite{Liu},
 and by myself when $L$ is fixed by an antisymplectic involution, see \cite{WelsCRAS}, \cite{WelsInvent}, \cite{WelsICM}
 or also \cite{Cho}, \cite{Solomon}. In my recent work \cite{Open4}, I define such open Gromov-Witten invariants when $X$ is four-dimensional
 and $L$ orientable. The invariance then only holds modulo $q \in \N$ for discs homologous to a class $d \in H_2 (X , L ; \Z)$ such that
  $\partial d = 0 \in H_1 (L ; \Z/q\Z)$. I also introduce similar invariants by counting reducible discs with a given number of irreducible components,
  see  \cite{Open4}. The aim of this paper is to define similar open Gromov-Witten invariants in a six-dimensional symplectic manifold. 
 
   Let thus $(X , \omega)$ be a  closed symplectic six-manifold and $L$ be a closed Lagrangian submanifold of $X$.
We again assume that $L$ is orientable and moreover that the inclusion of $L$ into $X$ induces an injective morphism
$H_1 (L ; A) \to H_1 (X ; A)$, where $A$ is a commutative ring. We also assume that $X$ contains no Maslov zero
pseudo-holomorphic disc with boundary on $L$, in order no to have to take into account branched cover of discs. Under these hypotheses, 
given a relative homology class $d \in H_2 (X , L ; \Z)$
of positive Maslov index, a generic almost-complex structure $J$ tamed by the symplectic form and submanifolds of $L$ and $X$ of adequate
cardinality and dimensions, we prove that the weighted number of $J$-holomorphic discs with boundary on $L$, homologous to $d$, either 
irreducible or reducible disconnected and which meet all of the chosen submanifolds of $L$ and $X$, only depends on the homology classes
of the submanifolds, of $d$ and of $\omega$ up to deformation, while it dos not depend on the generic choice of $J$. For this to be true, we
nevertheless assume that at least one submanifold has been chosen in $L$, which means that our discs at least contain one marked point on their
boundaries, see Theorem \ref{theoGW} and Corollary \ref{corGW}. These $J$-holomorphic discs are counted with respect to a sign, as usual
in the theory of Gromov-Witten invariants, but also with respect to a weight when they have more than one irreducible component. This weight is
defined in the following way. Let $D$ be a $J$-holomorphic disc with boundary on $L$ having $n>1$ disjoint  irreducible components. The boundary of
$D$ is an oriented link in $L$, each component of which has trivial homology class in $H_1 (L ; A)$. We then label every vertex of the complete graph
$K_n$ (having $n$ vertices) with one component of this link $\partial D$. Every edge of $K_n$ gets then decorated with the linking number of the
knots associated to its boundary vertices, since we have equipped $L$ with an orientation and even actually with a spin structure. 
These linking numbers take value in the ring $A$. Then, for every spanning subtree $T$ of $K_n$, we associate the product of the $n-1$
linking numbers associated to its $n-1$ edges and get a number $T_* \in A$. The sum of all these numbers $T_*$ over all spanning
subtrees of $K_n$ provides the weight under which we count the disc $D$. This weight does not depend on the labeling of the vertices of
$K_n$, we call it self-linking weight, see Definition \ref{defselflinkingweight}. As a consequence, the open Gromov-Witten invariants that
we define here take value in the commutative ring $A$. Recall that $K_n$ contains exactly $n^{n-2}$ spanning subtrees, a formula
established by J. J. Sylvester and A. Cayley, see \cite{Syl}, \cite{Cay}. 

The paper is organized as follows. In the first paragraph, we introduce the moduli spaces of pseudo-holomorphic discs and discuss some
of their properties, see in particular Proposition \ref{propboundaries}. The second paragraph is devoted to linking numbers, complete graphs
and spanning subtrees. We establish there Lemma \ref{lemmaforested}, a key property of self-linking weights. The last paragraph is devoted 
to the statement and proof of our result. \\

{\bf Acknowledgements:}

The research leading to these results has received funding from the European Community's Seventh Framework Progamme 
([FP7/2007-2013] [FP7/2007-2011]) under grant agreement $\text{n}\textsuperscript{o}$ [258204].
I am also grateful to F. Chapoton for pointing out the references \cite{Syl}, \cite{Cay} to me.

\section{Pseudo-holomorphic discs with boundary on a Lagrangian submanifold}
\subsection{The automorphism group of Poincar\'e's unit disc}
\label{subsectiondisc}

Let $\Delta = \{ z \in \C \; \vert \; \vert z \vert \leq 1 \}$ be the closed complex unit disc and $\H = \{ z \in \C \;  \vert \; \Im z  >0 \}$ be the upper half plane. 
Denote by $\stackrel{\circ}{\Delta}$ the interior of $\Delta$ and by $\overline{\H}$ the closure of $\H$ in the Riemann sphere. The homography
$z \in \overline{\H}\mapsto \frac{z-i}{z+i} \in \Delta$ is an isomorphism which we fix in order to identify $\Delta$ and $\overline{\H}$ in the sequel. 

\subsubsection{Orientation}

The group $\text{Aut} (\Delta)$ of automorphisms of $\Delta $ is isomorphic to $PSL_2 (\R) = \text{Aut} (\H)$. It acts transitively and without fixed point
on the product $\partial \Delta \times \stackrel{\circ}{\Delta}$ by $(\phi , z , \zeta ) \in \text{Aut} (\Delta) \times \partial \Delta \times \stackrel{\circ}{\Delta}
\mapsto (\phi (z) , \phi(\zeta) ) \in  \partial \Delta \times \stackrel{\circ}{\Delta}$. The product $ \partial \Delta \times \stackrel{\circ}{\Delta}$ being
canonically oriented, this action induces an orientation on   $\text{Aut} (\Delta)$ which we fix once for all. 

In the same way, $\text{Aut} (\Delta)$ acts transitively and without fixed point on the open subset of $(\partial \Delta)^3$ made of cyclically ordered
 triple of distinct points. The orientation induced on $\text{Aut} (\Delta)$ by this action coincides with the one we just fixed, see Lemma \ref{lemmaorient}. 
 
 For every $\epsilon >0$, denote by $r_\epsilon : z \in \Delta \mapsto \exp (i \epsilon) z  \in \Delta$, $t_\epsilon : z \in \H \mapsto z + \epsilon  \in \H$ and
 $h_\epsilon : z \in \H \mapsto (1 + \epsilon) z  \in \H$. Denote by $\stackrel{.}{r}_0 = \frac{\partial}{\partial \epsilon}\vert_{\epsilon = 0} r_\epsilon$,
 $\stackrel{.}{t}_0 = \frac{\partial}{\partial \epsilon}\vert_{\epsilon = 0} t_\epsilon$ and $\stackrel{.}{h}_0 = \frac{\partial}{\partial \epsilon}\vert_{\epsilon = 0} h_\epsilon$
 the associated holomorphic vector fields of  $\Delta$.
 
 \begin{lemma}
 \label{lemmaorient}
 The triple $(\stackrel{.}{r}_0 , \stackrel{.}{t}_0 , \stackrel{.}{h}_0)$ forms a direct basis of the Lie algebra $\text{aut} (\Delta)$. Moreover, the action of
  $\text{Aut} (\Delta)$ on  the open subset of $(\partial \Delta)^3$ made of cyclically ordered triple of distinct points preserves orientations.
 \end{lemma}
 
 {\bf Proof:}
 
 Let us identify $\partial \Delta \times \stackrel{\circ}{\Delta}$  with the orbit of $(1,0)$ under the action of $\text{Aut} (\Delta)$. The differential of this action
 maps $\stackrel{.}{r}_0$ onto the positive generator of $T_1  \partial \Delta$  whereas it maps $(\stackrel{.}{t}_0 , \stackrel{.}{h}_0)$ onto a direct basis
 of  $T_i \H$. The first part of the lemma follows. Now, $(\infty , -1 , 1)$ forms a cyclically ordered triple of distinct points on $\partial \overline{\H}$. The
 differential of the action of $ \text{Aut} (\H)$ on this triple maps $\stackrel{.}{r}_0 $ on the sum of positive generators of $T_\infty \partial  \overline{\H}$, $T_{-1} \partial \overline{\H}$ and
 $T_1 \partial  \overline{\H}$, whereas it maps $\stackrel{.}{t}_0$ onto the sum of positive generators of $T_{-1} \partial \overline{\H}$ and
 $T_1 \partial  \overline{\H}$ and $\stackrel{.}{h}_0$ onto the sum of a negative generator of $T_{-1} \partial \overline{\H}$ and a positive generator of
 $T_1 \partial  \overline{\H}$, hence the result. $\square$
 
 \begin{rem}
 It follows from Lemma \ref{lemmaorient} that our convention of orientation of $\text{Aut} (\Delta)$ is the same as the one adopted by Fukaya, Oh, Ohta and Ono in \S $8$ of \cite{FOOO}.
 \end{rem}
 
 \subsubsection{Glueing of automorphisms}
 \label{subsecglue}
 
 Denote by $\Delta_0$ the nodal disc obtained as the closure of $(\H \times \{ 0 \}) \cup (\{ 0 \}   \times \H ) \subset \C^2$ in $\C P^2$.
 Likewise, for every $\eta \geq 0$, denote by $\Delta_\eta$ the closure in $\C P^2$ of the hyperbola 
 $\{ (x,y) \in \H^2 \; \vert \; xy = -\eta \}$. For every $\epsilon >0$, denote by $\tau_\epsilon : x \in \H \mapsto \frac{x}{1-\epsilon x} \in \H$
 the translation of $\H$ fixing the origin  and by $\kappa_\epsilon : x \in \H \mapsto \frac{1}{1+ \epsilon } x \in \H$ the homothety of weight
 $ \frac{1}{1+ \epsilon }$, so that the corresponding pair $(\stackrel{.}{\tau}_0 , \stackrel{.}{\kappa}_0)$ of holomorphic vector fields
 forms a direct basis of the Lie algebra $\text{aut} (\overline{\H} , 0)$, where $\text{Aut} (\overline{\H} , 0)$ denotes the group of 
 automorphisms of $\overline{\H}$ fixing the boundary point $0$. The group of automorphisms of $\Delta_0$ is isomorphic to 
 $\text{Aut} (\overline{\H} , 0)^2$ and is oriented in such a way that its action on the interior $\stackrel{\circ}{\Delta}_0$ gets orientation preserving. 
 
 \begin{prop}
 \label{propnodal}
 The elements $(\stackrel{.}{\tau}_0 , 0)$ , $(0, \stackrel{.}{\tau}_0)$ and $(-\stackrel{.}{\kappa}_0 , \stackrel{.}{\kappa}_0)$ of the Lie algebra
 $\text{aut} (\Delta_0) = \text{aut} (\overline{\H} , 0)^2$ holomorphically deform as a direct basis of $\text{aut} (\Delta_\eta)$ for every $\eta >0$.
 Moreover, the infinitesimal action of $(\stackrel{.}{\kappa}_0 , 0)$ by reparametrization of maps $\Delta_0 \to \Delta_0$ extends as an inward
 normal holomorphic vector field transverse to the fibers $\Delta_\eta$ of $\cup_{\eta >0} \Delta_\eta$, lifting $\frac{\partial}{\partial \eta}$.
 \end{prop}

The action of $\text{Aut} (\Delta_0)$ by reparameterization of maps $\Delta_0 \to \Delta_0$ is defined by
$(g,u) \in \text{Aut} (\Delta_0) \times \text{Hol} (\Delta_0, \Delta_0) \mapsto u \circ g^{-1} \in \text{Hol} (\Delta_0, \Delta_0)$.\\

{\bf Proof:}

For every $\eta >0$, let us parameterize the interior of $\Delta_\eta$ by the map $x \in \H \mapsto (x , \frac{-\eta}{x}) \in \stackrel{\circ}{\Delta}_\eta$. 
For every $\epsilon \geq 0$, the translation $\tau_\epsilon $ induces through this parameterization the automorphism
$(x,y) \in \Delta_\eta \mapsto (\tau_\epsilon (x) , y + \epsilon \eta)  \in \Delta_\eta $. The associated holomorphic vector field $(\stackrel{.}{\tau}_0 , \eta)$
of $\Delta_\eta$ deforms $(\stackrel{.}{\tau}_0 , 0)$. By symmetry, the pair $(\eta , \stackrel{.}{\tau}_0)$ defines an element of $\text{aut} (\Delta_\eta)$
deforming  $(0 , \stackrel{.}{\tau}_0)$. Likewise, the homothety $\kappa_\epsilon$ induces through our parameterization the automorphism
$(x,y) \in \Delta_\eta \mapsto (\kappa_\epsilon (x) , h_\epsilon (y ))  \in \Delta_\eta $. The associated holomorphic vector field
$(\stackrel{.}{\kappa}_0 , -\stackrel{.}{\kappa}_0)$ deforms the opposite of $(- \stackrel{.}{\kappa}_0 , \stackrel{.}{\kappa}_0)$. The evaluation map at $i \in \H$
sends the pair $(\stackrel{.}{\tau}_0 , \stackrel{.}{\kappa}_0)$ to a direct basis of $T_i \H$. By deformation, for every $\eta >0$ close to zero, the pairs
$(\stackrel{.}{\tau}_0 , \eta)$ and $(\stackrel{.}{\kappa}_0 , -\stackrel{.}{\kappa}_0)$ evaluate as a direct basis of $T_{(i, \eta i)} \Delta_\eta $. The vector field
$(\eta , \stackrel{.}{\tau}_0)$ evaluates as a direct basis of $T \partial  \Delta_\eta $ and is close to zero at $(i, \eta i)$. As a consequence, the triple
$(\eta , \stackrel{.}{\tau}_0)$, $(\stackrel{.}{\tau}_0 , \eta)$ and $(\stackrel{.}{\kappa}_0 , -\stackrel{.}{\kappa}_0)$ defines a direct basis of $\text{aut} (\Delta_\eta)$.

 Finally, the biholomorphism $(x,y) \in \H^2 \mapsto (\eta x , y)  \in \H^2$ maps $\Delta_1$ onto $\Delta_\eta$ for every $\eta >0$. It is obtained by integration
 of the vector field $(\stackrel{.}{h}_0 , 0)$. The latter is the image of $(\stackrel{.}{\kappa}_0 , 0)$ under the infinitesimal action of $\text{Aut} (\Delta_0)$
 by reparameterization of maps $\Delta_0 \to \Delta_0$. $\square$
 
 \begin{cor}
 \label{cornodal}
 Let $( \Delta_\eta)_{\eta \geq 0}$ be the standard deformation of the nodal disc $\Delta_0$. Then , a direct basis of the Lie algebra $\text{aut} (\Delta_0)$
 deforms as the concatenation of an outward normal vector field lifting $- \frac{\partial}{\partial \eta}$ and a direct basis of $\text{aut} (\Delta_\eta)$, $\eta >0$.
 \end{cor}

{\bf Proof:}

This Corollary \ref{cornodal} follows from Proposition \ref{propnodal} and the fact that the quadruple
$(\stackrel{.}{\tau}_0 , 0)$ , $(0, \stackrel{.}{\tau}_0)$, $(-\stackrel{.}{\kappa}_0 , \stackrel{.}{\kappa}_0)$ and $(\stackrel{.}{\kappa}_0 , 0)$
defines a direct basis of $\text{aut} (\Delta_0)$. $\square$

\subsection{Moduli spaces of simple discs}
\label{subsecsimple}

We denote by
${\cal P} (X,L) = \{ (u,J) \in C^{1} (\Delta , X) \times {\cal J}_\omega \; \vert \; u(\partial \Delta) \subset L \text{ and } du + J\vert_u \circ du \circ j_{st} = 0\}$
the space of pseudo-holomorphic maps from $\Delta$ to the pair $(X,L)$, where $ j_{st}$ denotes the standard complex structure of $\Delta$
and ${\cal J}_\omega$ the space of almost-complex structures on $X$ of class $C^l$ tamed by $\omega$, $l \gg 1$.
Note that $J$ being of class $C^l$,  the regularity of such pseudo-holomorphic maps $u$ is actually more than $C^l$, see \cite{McDSal}. 
More generally, for every $r,s \in \N$, we denote by ${\cal P}_{r,s} (X,L) = \{ ((u,J), \underline{z}, \underline{\zeta}) \in {\cal P} (X,L) \times ((\partial \Delta)^r \setminus
\text{diag}_{\partial \Delta}) \times ((\stackrel{\circ}{\Delta})^s \setminus
\text{diag}_{\Delta}) \}$, where $\text{diag}_{\partial \Delta} = \{ (z_1, \dots , z_r) \in (\partial \Delta)^r \; \vert \; \exists i \neq j , z_i = z_j \}$ and
$\text{diag}_{ \Delta} = \{ (\zeta_1, \dots , \zeta_s) \in \Delta^s \; \vert \; \exists i \neq j , \zeta_i = \zeta_j \}$.

Following \cite{Laz}, \cite{KwonOh}, \cite{BirCor}, we define
\begin{defi}
\label{defsimple}
A pseudo-holomorphic map $u$ is said to be simple iff there is a dense open subset $\Delta_{inj} \subset \Delta$ such that
$\forall z \in \Delta_{inj}, u^{-1} (u(z)) = \{ z \}$ and $du\vert_z \neq 0$.
\end{defi}

We denote by ${\cal P}_{r,s}^* (X,L) $ the subset of simple elements of ${\cal P}_{r,s} (X,L) $. It is a separable Banach manifold which is
naturally embedded as a submanifold of class $C^{l-k}$ of the space $W^{k,p} (\Delta , X)  \times {\cal J}_\omega $ for every
$1 \ll k \ll l$ and $p > 2$, see Proposition $3.2$ of \cite{McDSal}. 

For every $d \in H_2 (X , L ; \Z)$, we denote by ${\cal P}_{r,s}^d (X,L) = \{ (u,J) \in {\cal P}_{r,s}^* (X,L)  \; \vert \; u_* [\Delta ] = d \}$ and by
${\cal M}_{r,s}^d (X,L)  = {\cal P}_{r,s}^d (X,L) /\text{Aut}(\Delta) $, where $\text{Aut}(\Delta)$ acts by composition on the right, see \S \ref{subsectiondisc}. 
The latter is equipped with a projection $\pi : [u,J, \underline{z}, \underline{\zeta}] \in {\cal M}_{r,s}^d (X,L) \mapsto
J \in {\cal J}_\omega$ and an evaluation map $\text{ev} :  [u,J, \underline{z}, \underline{\zeta}] \in {\cal M}_{r,s}^d (X,L) \mapsto (u(\underline{z}), u(\underline{\zeta}))
\in L^r \times X^s$. 

We recall the following classical result due to Gromov (see \cite{Gro}, \cite{McDSal}, \cite{FOOO}).
\begin{theo}
\label{theodim}
For every closed Lagrangian submanifold $L$ of a six-dimensional closed symplectic manifold $(X , \omega)$
and every $d \in H_2 (X , L ; \Z)$, $r,s \in \N$, the space ${\cal M}_{r,s}^d (X,L)$ is a separable Banach manifold and the projection
$\pi :  {\cal M}_{r,s}^d (X,L) \to {\cal J}_\omega$ is Fredholm of index $\mu_L (d) + r +2s$. $\square$
\end{theo}
Note that from Sard-Smale's theorem \cite{Smale}, the set of regular values of $\pi$ is dense of the second category. As a consequence,
for a generic choice of $J \in  {\cal J}_\omega$, the moduli space ${\cal M}_{0,0}^d (X,L ; J) = \pi^{-1} (J)$ is a manifold of dimension
$\mu_L (d)$ as soon as it is not empty. 

Forgetting the marked points defines a map $(u,J, \underline{z}, \underline{\zeta}) \in {\cal P}_{r,s}^d (X,L) \mapsto (u,J) \in {\cal P}_{0,0}^d (X,L)$
which quotients as a forgetful map $f_{r,s} : {\cal M}_{r,s}^d (X,L) \to {\cal M}_{0,0}^d (X,L) $ whose fibers are canonically oriented. 
When $L$ is equipped with a Spin structure, the manifolds ${\cal M}_{r,s}^d (X,L ; J) $, $r,s \in \N$,  inherit canonical orientations, see Theorem $8.1.1$ of \cite{FOOO}. 
We adopt the following convention for orienting these moduli spaces. 
When $J$ is a generic almost-complex structure tamed by $\omega$, the manifold $ {\cal P}_{0,0}^d (X,L ; J)$ inherits an orientation from the Spin
structure of $L$, see Theorem $8.1.1$ of \cite{FOOO}. This orientation induces an orientation on the manifold 
${\cal M}_{0,0}^d (X,L ; J) = {\cal P}_{0,0}^d (X,L ; J) /\text{Aut}(\Delta) $
such that at every point $(u,J) \in {\cal P}_{0,0}^d (X,L ; J)$, the concatenation of a direct basis of $T_{[u,J]} {\cal M}_{0,0}^d (X,L ; J)$ in a horizontal space of
$T_{(u,J)} {\cal P}_{0,0}^d (X,L ; J)$ followed by a direct basis of the orbit of $\text{Aut}(\Delta)$ at $(u,J)$ defines a direct basis of $T_{(u,J)} {\cal P}_{0,0}^d (X,L ; J)$.
This is the convention $8.2.1.2$ adopted by Fukaya, Oh, Ohta and Ono in
 \cite{FOOO}, so that our orientation of ${\cal M}_{0,0}^d (X,L ; J)$, which we call the quotient orientation, coincides with the one introduced in
 \cite{FOOO}. The manifold ${\cal M}_{r,s}^d (X,L ; J) $ is then oriented in such a way that at every point, the concatenation of a direct basis of $T_{[u,J]} {\cal M}_{0,0}^d (X,L ; J)$
 with a direct basis of the fiber of the forgetful map $f_{r,s}$ provides a direct basis of $T_{[u,J]} {\cal M}_{r,s}^d (X,L ; J)$. 
 (The latter convention differs in general from the one adopted by Fukaya, Oh, Ohta and Ono in  \cite{FOOO}.)

\subsection{Moduli spaces of nodal discs}
\label{subsecred}

Let $d_1 , d_2 \in H_2 (X, L ; \Z)$ be such that $\mu_L (d_1) >1$, $\mu_L (d_2) >1$. For every $z \in \partial \Delta$ and $i \in \{1,2\}$, denote by 
$ev_z : (u,J) \in {\cal P}_{0,0}^{d_i} (X,L) \mapsto u(z) \in L$ the evaluation map at the point $z$. As soon as $J$ is generic, the restrictions
of these evaluation maps to $ {\cal P}_{0,0}^{d_i} (X,L ; J)$ are submersions.  We then denote by
$ {\cal P}_{0,0}^{(d_1 , d_2)} (X,L)$ the fiber product ${\cal P}_{0,0}^{d_1} (X,L) \text{    }_{ev_{1}} \! \! \! \times_{ev_{-1}} {\cal P}_{0,0}^{d_2} (X,L)$.
As soon as $J$ is generic enough, $ {\cal P}_{0,0}^{(d_1 , d_2)} (X,L ; J)$ is a manifold of dimension $\mu_L (d_1 + d_2) + 3$ since $X$ is six-dimentional
throughout this paper. When $L$ is equipped with a Spin structure, it is canonically oriented from \cite{FOOO}. Let us recall the convention of orientation of this manifold. 

Let $(u_1 , u_2 , J) \in {\cal P}_{0,0}^{(d_1 , d_2)} (X,L ; J)$ and $(v_1 , v_2 , v_3)$ (resp. $(w_1 , w_2 , w_3)$) be elements of 
$T_{u_1} {\cal P}_{0,0}^{d_1} (X,L ; J) $ (resp. $T_{u_2} {\cal P}_{0,0}^{d_2} (X,L ; J) $) such that $(v_1 (1) , v_2 (1) , v_3 (1))$ (resp. $(w_1 (-1) , w_2  (-1), \\
 w_3  (-1))$)  forms a direct basis of $T_{u_1 (1)} L$ (resp. $w_1 (-1) = v_1 (1)$, $w_2 (-1) = v_2 (1)$, $w_3 (-1) = v_3 (1)$).
Let ${\cal B}_1$ (resp.  ${\cal B}_2$) be an ordered family of $\mu_L (d_1) $ (resp. $\mu_L (d_2) $) elements of 
$T_{u_1} {\cal P}_{0,0}^{d_1} (X,L ; J) $ (resp. $T_{u_2} {\cal P}_{0,0}^{d_2} (X,L ; J) $) such that $({\cal B}_1 , v_1 , v_2 , v_3)$ (resp. $(w_1 , w_2 , w_3 , {\cal B}_2)$)
defines a direct basis of $T_{u_1} {\cal P}_{0,0}^{d_1} (X,L ; J) $ (resp. $T_{u_2} {\cal P}_{0,0}^{d_2} (X,L ; J) $).
Then, the manifold $ {\cal P}_{0,0}^{(d_1 , d_2)} (X,L ; J)$ is oriented such  that the basis $({\cal B}_1 , v_1 + w_1 , v_2 + w_2 , v_3 + w_3 , {\cal B}_2)$ of  \\
$T_{(u_1 , u_2 , J)} {\cal P}_{0,0}^{(d_1 , d_2)} (X,L ; J)$ becomes direct. We then denote by $ {\cal M}_{0,0}^{(d_1 , d_2)} (X,L)$ the quotient
$ {\cal P}_{0,0}^{(d_1 , d_2)} (X,L) / \text{Aut} (\Delta_0)$, where $\Delta_0$ denotes the nodal disc, see \S \ref{subsecglue}.
When $J$ is generic enough, $ {\cal M}_{0,0}^{(d_1 , d_2)} (X,L)$ is a manifold of dimension  $\mu_L (d_1 + d_2) - 1$ which we equip with the quotient orientation
as in \S \ref{subsecsimple}. This convention of orientation coincides thus with the one adopted in \cite{FOOO}. Note the the tautological involution
$ {\cal M}_{0,0}^{(d_1 , d_2)} (X,L) \to  {\cal M}_{0,0}^{(d_2 , d_1)} (X,L)$ which exchanges the discs preserves the orientation since the Maslov indices
$\mu_L (d_1)$, $\mu_L (d_2)$ are even. Finally, for every $r,s \in \N$, we denote by ${\cal P}_{r,s}^{(d_1 , d_2)} (X,L) = {\cal P}_{0,0}^{(d_1 , d_2)} (X,L)
\times  ((\partial \Delta_0  \setminus \{ \text{node} \})^r \setminus
\text{diag}_{\partial \Delta_0}) \times ((\stackrel{\circ}{\Delta}_0)^s \setminus
\text{diag}_{\Delta_0}) \}$, where $\text{diag}_{\partial \Delta_0} = \{ (z_1, \dots , z_r) \in (\partial \Delta_0 )^r \; \vert \; \exists i \neq j , z_i = z_j \}$ and
$\text{diag}_{ \Delta_0} = \{ (\zeta_1, \dots , \zeta_s) \in \Delta_0^s \; \vert \; \exists i \neq j , \zeta_i = \zeta_j \}$. We then denote by 
${\cal M}_{r,s}^{(d_1 , d_2)} (X,L) = {\cal P}_{r,s}^{(d_1 , d_2)} (X,L) / \text{Aut} (\Delta_0)$ and by $f_{r,s} : {\cal M}_{r,s}^{(d_1 , d_2)} (X,L) \to {\cal M}_{0,0}^{(d_1 , d_2)} (X,L) $
the associated forgetful map, whose fibers are canonically oriented. 
When $J$ is a generic almost-complex structure tamed by $\omega$,  the manifolds ${\cal M}_{r,s}^{(d_1 , d_2)} (X,L ; J) $ are oriented in such a way that at every point
$[u_1 , u_2 ,J , \underline{z}, \underline{\zeta}] \in {\cal M}_{r,s}^{(d_1 , d_2)} (X,L ; J) $,
the concatenation of a direct basis of $T_{[u_1 , u_2 ,J , \underline{z}, \underline{\zeta}]} {\cal M}_{0,0}^{(d_1 , d_2)} (X,L ; J)$ in a horizontal space 
 with a direct basis of the fiber of the forgetful map $f_{r,s}$ at $[u_1 , u_2 ,J, \underline{z}, \underline{\zeta}]$ provides a direct basis of 
 $T_{[u_1 , u_2 ,J, \underline{z}, \underline{\zeta}]} {\cal M}_{r,s}^{(d_1 , d_2)} (X,L ; J)$. (The latter convention differs in general from
 the one adopted by Fukaya, Oh, Ohta and Ono in  \cite{FOOO}.)
 
 From Gromov's compactness and glueing theorems (see for example \cite{Frauen}, \cite{FOOO}, \cite{BirCor}), the space ${\cal M}_{r,s}^{(d_1 , d_2)} (X,L ; J) $ canonically identifies  as a component of the boundary
 of the moduli space ${\cal M}_{r,s}^{d_1 + d_2} (X,L ; J) $. The following Proposition \ref{propboundary1}, analogous to Proposition $8.3.3$ of \cite{FOOO},
  compares the orientation of ${\cal M}_{r,s}^{(d_1 , d_2)} (X,L ; J) $ with the one induced by ${\cal M}_{r,s}^{d_1 + d_2} (X,L ; J) $.
 
 \begin{prop}
 \label{propboundary1}
 Let $L$ be a closed Lagrangian Spin submanifold of a closed symplectic six-manifold $(X , \omega)$.
 Let $d_1 , d_2 \in H_2 (X, L ; \Z)$ be such that $\mu_L (d_1) \geq 2$, $\mu_L (d_2) \geq 2$ and let $r,s \in \N$. 
 Then, for every generic almost-complex structure $J$ tamed by $\omega$, the incidence index
 $\langle \partial {\cal M}_{r,s}^{d_1 + d_2} (X,L ; J)  , {\cal M}_{r,s}^{(d_1 , d_2)} (X,L ; J)  \rangle$ equals $-1$. 
 \end{prop}

{\bf Proof:}

The glueing map of $J$-holomorphic discs preserves the orientations of the fibers of the forgetful map $f_{r,s}$, so that from our conventions
of orientations of moduli spaces, it suffices to prove the result for $r=s=0$. From Lemma $8.3.5$ of  \cite{FOOO},
the glueing map of $J$-holomorphic maps  ${\cal P}_{0,0}^{(d_1 , d_2)} (X,L ; J) \to {\cal P}_{0,0}^{d_1 + d_2} (X,L ; J) $ preserves orientations. 
Let $(u_1 , u_2 ,J)  \in {\cal P}_{0,0}^{(d_1 , d_2)} (X,L ; J) $, ${\cal B}_1$ be a direct basis of $T_{[u_1 , u_2 ,J]} {\cal M}_{0,0}^{(d_1 , d_2)} (X,L ; J)$ 
and ${\cal B}_0$ a direct basis of the linearized orbit of $ \text{Aut} (\Delta_0)$ at $(u_1 , u_2 ,J) $. Then, by definition, the concatenation
$({\cal B}_1 , {\cal B}_0)$ defines a direct basis of $T_{(u_1 , u_2 ,J)}  {\cal P}_{0,0}^{(d_1 , d_2)} (X,L ; J) $. Now, let $(u_1 \#_R u_2)_{R \gg 1}$ be
a path of ${\cal P}_{0,0}^{d_1 + d_2} (X,L ; J) $ transversal to ${\cal P}_{0,0}^{(d_1 , d_2)} (X,L ; J) $ at $u_1 \#_\infty u_2 = (u_1 , u_2)$
given by the glueing map. From Corollary \ref{cornodal}, ${\cal B}_0$ deforms for $R < + \infty$ as a pair $(n , {\cal B}'_0)$ of
$T_{(u_1 \#_R u_2 , J)} {\cal P}_{0,0}^{d_1 + d_2} (X,L ; J) $, where ${\cal B}'_0$ is a direct basis of the  linearized orbit of $ \text{Aut} (\Delta)$ at
$u_1 \#_R u_2$ and $n$ points towards the boundary of ${\cal P}_{0,0}^{d_1 + d_2} (X,L ; J) $. As a consequence,
$({\cal B}_1 , n)$ deforms as a direct basis of $T_{(u_1 \#_R u_2 , J)} {\cal M}_{0,0}^{d_1 + d_2} (X,L ; J) $. The result now follows from the fact that
the cardinality of ${\cal B}_1$ is odd. $\square$

\subsection{Moduli spaces of reducible discs}
\label{subsecreducible}

Let $n \geq 1$ and $d_1 , \dots , d_n \in H_2 (X, L ; \Z)$ be such that $\mu_L (d_i) >0$, $i \in \{ 1, \dots , n \}$. For 
every almost-complex structure $J$ tamed by $\omega$, we denote by ${\cal M}_{0,0}^{d_1, \dots , d_n} (X,L ; J) $
the direct product ${\cal M}_{0,0}^{d_1} (X,L ; J) \times \dots \times {\cal M}_{0,0}^{d_n} (X,L ; J) $ and by 
${\cal P}_{0,0}^{d_1, \dots , d_n} (X,L) $ the corresponding fiber product ${\cal P}_{0,0}^{d_1} (X,L) \times_{{\cal J}_\omega} \dots \times_{{\cal J}_\omega}
 {\cal P}_{0,0}^{d_n} (X,L) $. When $J$ is generic enough and $L$ Spin, ${\cal M}_{0,0}^{d_1, \dots , d_n} (X,L ; J) $ is a manifold of dimension
 $\mu_L (d_1 + \dots + d_n)$ equipped with its product orientation. Since the manifolds ${\cal M}_{0,0}^{d_i} (X,L ; J)$ are even dimensional,
 $i \in \{ 1 , \dots , n \}$, this product orientation coincides with the quotient orientation induced by ${\cal P}_{0,0}^{d_1, \dots , d_n} (X,L ; J) /  \text{Aut}(\Delta)^n$,
 where ${\cal P}_{0,0}^{d_1, \dots , d_n} (X,L ; J) = \prod_{i=1}^n {\cal P}_{0,0}^{d_i} (X,L ; J)$ is itself equipped with the product orientation. 
For every $r,s \in \N$, we likewise define ${\cal M}_{r,s}^{d_1, \dots , d_n} (X,L) $ as the quotient of
${\cal P}_{0,0}^{d_1, \dots , d_n} (X,L)  \times ((\partial \Delta \cup \dots \cup \partial \Delta)^r \setminus
\text{diag}_{\partial \Delta}) \times ((\stackrel{\circ}{\Delta}  \cup \dots \cup \stackrel{\circ}{\Delta})^s \setminus
\text{diag}_{\Delta}) \}$ by $ \text{Aut}(\Delta)^n$, where $\text{diag}_{\partial \Delta} = \{ (z_1, \dots , z_r) \in (\partial \Delta  \cup \dots \cup \partial \Delta)^r \; \vert \; 
\exists i \neq j , z_i = z_j \}$ and
$\text{diag}_{ \Delta} = \{ (\zeta_1, \dots , \zeta_s) \in (\Delta  \cup \dots \cup \Delta)^s \; \vert \; \exists i \neq j , \zeta_i = \zeta_j \}$. 
It is equipped with a forgetful map $f_{r,s} : {\cal M}_{r,s}^{d_1, \dots , d_n} (X,L)  \to {\cal M}_{0,0}^{d_1, \dots , d_n} (X,L) $ whose fibers
are canonically oriented by the complex structure. For generic $J$, the manifold ${\cal M}_{r,s}^{d_1, \dots , d_n} (X,L ; J) $ is as before oriented
in such a way that at every point $[u_1 , \dots , u_n, J, \underline{z}, \underline{\zeta}]$, 
the concatenation of a direct basis of $T_{[u_1 , \dots , u_n,J, \underline{z}, \underline{\zeta}]} {\cal M}_{0,0}^{d_1, \dots , d_n} (X,L ; J) $
in a horizontal space with a direct basis of the fiber of the forgetful map $f_{r,s}$ at $[u_1 , \dots , u_n,J, \underline{z}, \underline{\zeta}]$ 
provides a direct basis of $T_{[u_1 , \dots , u_n,J]} {\cal M}_{r,s}^{d_1, \dots , d_n} (X,L ; J)$. This
 orientation is not the quotient orientation of ${\cal P}_{r,s}^{d_1, \dots , d_n} (X,L ; J) /  \text{Aut}(\Delta)^n$.
 The tautological action of the group of permutation of $\{Êd_1, \dots , d_n \}$ on these spaces ${\cal M}_{r,s}^{d_1, \dots , d_n} (X,L ; J) $ preserves the orientations,
 since the manifolds ${\cal M}_{0,0}^{d_i} (X,L ; J)$ are even dimensional.
 
 In the same way, for every $n \geq 2$, we denote by ${\cal M}_{0,0}^{(d_1, d_2), d_3, \dots , d_n} (X,L ; J) $
the direct product ${\cal M}_{0,0}^{(d_1, d_2)} (X,L ; J) \times {\cal M}_{0,0}^{d_3} (X,L ; J) \times \dots \times {\cal M}_{0,0}^{d_n} (X,L ; J) $
and equip it with the product orientation, see \S \ref{subsecred}. We likewise define, for every $r,s \in \N$, the corresponding space 
${\cal M}_{r,s}^{(d_1, d_2), d_3, \dots , d_n} (X,L ; J) $ and still denote by $f_{r,s} : {\cal M}_{r,s}^{(d_1, d_2), d_3, \dots , d_n} (X,L ; J) 
\to {\cal M}_{0,0}^{(d_1, d_2), d_3, \dots , d_n} (X,L ; J) $ the forgetful map, whose fibers are canonically oriented. We orient 
${\cal M}_{r,s}^{(d_1, d_2), d_3, \dots , d_n} (X,L ; J) $ as the concatenation of the orientation of ${\cal M}_{0,0}^{(d_1, d_2), d_3, \dots , d_n} (X,L ; J) $ 
with the orientation of the fibers of $f_{r,s} $.

Note that every element of ${\cal M}_{0,0}^{(d_1, d_2)} (X,L ; J) $ is parameterized by the nodal disc $\Delta_0$ and thus possess a
special point, the node $\bullet$ of  $\Delta_0$. We denote by $f_\bullet : {\cal M}_{0,0}^{(d_1, d_2)} (X,L ; J) \to {\cal M}_{0,0}^{d_1, d_2} (X,L ; J) $
the tautological map forgetting this special point $\bullet$ and also by $f_\bullet : {\cal M}_{r,s}^{(d_1, d_2), d_3, \dots , d_n} (X,L ; J) 
\to {\cal M}_{r,s}^{d_1, \dots , d_n} (X,L ; J) $ the induced forgetful map. 

\begin{lemma}
\label{lemmawall}
Let $L$ be a closed oriented Lagrangian submanifold of a closed symplectic six-manifold $(X , \omega)$.
 Let $r,s \in \N$,   $n \geq 2$ and $d_1, \dots , d_n \in H_2 (X, L ; \Z)$ be such that $\mu_L (d_j) \geq 2$, $j \in \{1, \dots , n \}$. 
 Then, for every generic almost-complex structure $J$ tamed by $\omega$, the image $f_\bullet \big( {\cal M}_{r,s}^{(d_1, d_2), d_3, \dots , d_n} (X,L ; J) \big)$
 is, outside of a codimension two subspace, a canonically cooriented codimension one submanifold of $ {\cal M}_{r,s}^{d_1, \dots , d_n} (X,L ; J) $.
 \end{lemma}

Note that the Lagrangian $L$ is only supposed to be oriented in Lemma \ref{lemmawall}, so that the canonical coorientation given by
this lemma does not originate from orientations of the moduli spaces ${\cal M}_{r,s}^{(d_1, d_2), d_3, \dots , d_n} (X,L ; J)$
and  $ {\cal M}_{r,s}^{d_1, \dots , d_n} (X,L ; J) $.\\

{\bf Proof:}

Let $[u_1 , \dots , u_n, J , \underline{z} , \underline{\zeta}] \in {\cal M}_{r,s}^{(d_1, d_2), d_3, \dots , d_n} (X,L ; J)$ where $J$ is generic.
Outside of a subspace of codimension more than two of $ {\cal M}_{r,s}^{(d_1, d_2), d_3, \dots , d_n} (X,L ; J)$, the differentials
of $u_1$ and $u_2$ at the node $\bullet $ do not vanish. Moreover, outside of a codimension two subspace of 
 $ {\cal M}_{r,s}^{(d_1, d_2), d_3, \dots , d_n} (X,L ; J)$, the lines $\text{Im} (d\vert_{\bullet} u_1\vert_{\partial \Delta})$ and
 $\text{Im} (d\vert_{\bullet} u_2 \vert_{\partial \Delta})$ are in direct sum in $T_{u_1 (\bullet)} L$. These lines are moreover canonically oriented,
 so that the normal line $N_\bullet$ to the oriented plane $\text{Im} (d\vert_{\bullet} u_1\vert_{\partial \Delta}) \oplus 
 \text{Im} (d\vert_{\bullet} u_2 \vert_{\partial \Delta})$ in the oriented three-space $T_{u_1 (\bullet)} L$ inherits an orientation. 
 Now, since the Maslov index $\mu_L (d_1)$ is greater than one and $J$ is generic, the map $u_1$ deforms as a family
 $(u_1^\lambda)_{\lambda \in ]-\epsilon, \epsilon[}$ in ${\cal M}_{0,0}^{d_1} (X,L ; J)$ such that $u_1^\lambda (\bullet)$
 is positively transverse to $\text{Im} (d\vert_{\bullet} u_1\vert_{\partial \Delta}) \oplus 
 \text{Im} (d\vert_{\bullet} u_2 \vert_{\partial \Delta})$ at $\lambda = 0$. The family 
 $[u_1^\lambda , u_2 , \dots , u_n, J , \underline{z} , \underline{\zeta}]_{\lambda \in ]-\epsilon, \epsilon[}$ of $ {\cal M}_{r,s}^{d_1, \dots , d_n} (X,L ; J)$ 
 is then transversal to $f_\bullet : {\cal M}_{r,s}^{(d_1, d_2), d_3, \dots , d_n} (X,L ; J) $ at $\lambda = 0$ and the canonical coorientation
 of this image is defined such that this family becomes positively transverse. $\square$\\

Note that the tautological involution which exchange the discs of class $d_1$ and $d_2$ preserves the coorientations given by
Lemma \ref{lemmawall}. When $L$ is Spin, the orientations of the spaces ${\cal M}_{r,s}^{(d_1, d_2), d_3, \dots , d_n} (X,L ; J)$
and $ {\cal M}_{r,s}^{d_1, \dots , d_n} (X,L ; J)$  also induce a coorientation of $f_\bullet : {\cal M}_{r,s}^{(d_1, d_2), d_3, \dots , d_n} (X,L ; J) $
in $ {\cal M}_{r,s}^{d_1, \dots , d_n} (X,L ; J)$. We denote the incidence index between these coorientations by 
$\langle {\cal M}_{r,s}^{d_1, \dots , d_n} (X,L ; J) ,  {\cal M}_{r,s}^{(d_1, d_2), d_3, \dots , d_n} (X,L ; J) \rangle$, so that it equals $+1$ 
if they coincide and $-1$ otherwise. 

\begin{prop}
\label{propboundaries}
Let $L$ be a closed Spin Lagrangian submanifold of a closed symplectic six-manifold $(X , \omega)$.
 Let $r,s \in \N$,   $n \geq 2$ and $d_1, \dots , d_n \in H_2 (X, L ; \Z)$ be such that $\mu_L (d_j) \geq 2$, $j \in \{1, \dots , n \}$. 
 Then, for every generic almost-complex structure $J$ tamed by $\omega$,
 $\langle {\cal M}_{r,s}^{d_1, \dots , d_n} (X,L ; J) ,  {\cal M}_{r,s}^{(d_1, d_2), d_3, \dots , d_n} 
(X,L ; J) \rangle = +1$, whereas
$$\langle \partial {\cal M}_{r,s}^{d_1 + d_2, d_3, \dots , d_n} (X,L ; J)  , {\cal M}_{r,s}^{(d_1 , d_2), d_3, \dots , d_n} (X,L ; J)  \rangle = -1.$$ 
\end{prop}

{\bf Proof:}

The glueing map of $J$-holomorphic discs preserves the orientations of the fibers of the forgetful map $f_{r,s}$, so that from our conventions
of orientations of moduli spaces, it suffices to prove the result for $r=s=0$. Moreover, the spaces ${\cal M}_{0,0}^{d_j} (X,L ; J)$ being even
dimensional, it suffices to prove the result for $n=2$. The second part of Proposition \ref{propboundaries} thus follows from Proposition
\ref{propboundary1}. Let then $[u_1 , u_2 ,J]  \in {\cal M}_{0,0}^{(d_1 , d_2)} (X,L ; J) $. We may assume that the lines 
$\text{Im} (d\vert_{\bullet} u_1\vert_{\partial \Delta})$ and
 $\text{Im} (d\vert_{\bullet} u_2 \vert_{\partial \Delta})$ are in direct sum in $T_{u_1 (\bullet)} L$ and that the differentials of the evaluation
 maps $ {\cal M}_{0,0}^{d_j} (X,L ; J) \to L$, $j \in \{ 1,2 \}$, at the node $\bullet$ are surjective. Let $\stackrel{.}{r}_1$ (resp.
 $\stackrel{.}{r}_2$) be an element of $\text{Aut} (\Delta)$ whose action on  ${\cal M}_{0,0}^{d_1} (X,L ; J)$ (resp.  ${\cal M}_{0,0}^{d_2} (X,L ; J)$)
 evaluated at the node $\bullet$ positively generates $\text{Im} (d\vert_{\bullet} u_1\vert_{\partial \Delta})$  (resp.  
 $\text{Im} (d\vert_{\bullet} u_2 \vert_{\partial \Delta})$). Let $r_1^*$ (resp. $r_2^*$) be an element of $T_{u_2} {\cal M}_{0,0}^{d_2} (X,L ; J)$
 (resp. $T_{u_1} {\cal M}_{0,0}^{d_1} (X,L ; J)$) whose evaluation at the node $\bullet$ coincides with $\stackrel{.}{r}_1$ (resp.
 $\stackrel{.}{r}_2$), so that $\stackrel{.}{r}_1 + r_1^*$ and $r_2^* + \stackrel{.}{r}_2$ define elements of 
 $T_{[u_1 , u_2 ,J] } {\cal M}_{0,0}^{(d_1 , d_2)} (X,L ; J) $. Let $\nu = \nu_1 + \nu_2 \in  {\cal M}_{0,0}^{(d_1 , d_2)} (X,L ; J) $ be an element such that $\nu (\bullet) = \nu_1(\bullet) = \nu_2 (\bullet)$
 positively generates the normal $N_\bullet$ to the oriented plane $\text{Im} (d\vert_{\bullet} u_1\vert_{\partial \Delta}) \oplus 
 \text{Im} (d\vert_{\bullet} u_2 \vert_{\partial \Delta})$ in the oriented three-space $T_{u_1 (\bullet)} L$. By construction, the basis
 $(\stackrel{.}{r}_1 \!\!\!(\bullet) , r_2^* (\bullet) ,\nu_1(\bullet) ) = ( r_1^* (\bullet) , \stackrel{.}{r}_2  \!\!\!(\bullet) , \nu_2 (\bullet) )$ of $T_{u_1 (\bullet)} L$
 is direct. Let $(\stackrel{.}{t}_1 , \stackrel{.}{h}_1)$ (resp. $(\stackrel{.}{t}_2 , \stackrel{.}{h}_2)$) be a pair of elements of $\text{aut} (\Delta)$ 
 which vanish at the node $\bullet$, so that $(\stackrel{.}{r}_1 , \stackrel{.}{t}_1 , \stackrel{.}{h}_1)$ (resp. $(\stackrel{.}{r}_2 , 
 \stackrel{.}{t}_2 , \stackrel{.}{h}_2)$) is a direct basis of  the Lie algebra $\text{aut} (\Delta)$. Let ${\cal B}_1$ (resp. ${\cal B}_2$) be a family of elements
 of $T_{(u_1 , J)} {\cal P}_{0,0}^{d_1} (X,L ; J)$ (resp. $T_{(u_2 , J)} {\cal P}_{0,0}^{d_2} (X,L ; J)$) such that the basis
 $({\cal B}_1 , \stackrel{.}{t}_1 , \stackrel{.}{h}_1 , \stackrel{.}{r}_1 , r_2^* ,\nu_1)$ (resp. $(r_1^*  , \stackrel{.}{r}_2 , \nu_2 , {\cal B}_2 , 
 \stackrel{.}{t}_2 , \stackrel{.}{h}_2)$) of $T_{(u_1 , J)} {\cal P}_{0,0}^{d_1} (X,L ; J)$ (resp. $T_{(u_2 , J)} {\cal P}_{0,0}^{d_2} (X,L ; J)$)  is direct. 
 By definition, the basis $({\cal B}_1 , \stackrel{.}{t}_1 , \stackrel{.}{h}_1 , \stackrel{.}{r}_1 + r_1^* , r_2^* + \stackrel{.}{r}_2 ,  \nu_1 + \nu_2 ,
 {\cal B}_2 , \stackrel{.}{t}_2 , \stackrel{.}{h}_2)$ of $T_{(u_1, u_2 , J)} {\cal P}_{0,0}^{(d_1 , d_2)} (X,L ; J)$ is direct, so that
 $({\cal B}_1 , \stackrel{.}{r}_1 + r_1^* , r_2^* + \stackrel{.}{r}_2 ,  \nu_1 + \nu_2 , {\cal B}_2)$ is a direct basis of 
 $T_{[u_1, u_2 , J]} {\cal M}_{0,0}^{(d_1 , d_2)} (X,L ; J)$. Likewise, $({\cal B}_1 , \stackrel{.}{t}_1 , \stackrel{.}{h}_1 , \stackrel{.}{r}_1 , r_2^* ,  \nu_1 ,
  r_1^* ,  \stackrel{.}{r}_2 ,  \nu_2 ,  {\cal B}_2 , \stackrel{.}{t}_2 , \stackrel{.}{h}_2)$ is a direct basis of $T_{(u_1, u_2 , J)} {\cal P}_{0,0}^{d_1 , d_2} (X,L ; J)$,
  so that $({\cal B}_1 ,  r_2^* ,  \nu_1 , r_1^* ,  \nu_2 ,  {\cal B}_2 )$ is an indirect basis of  $T_{[u_1, u_2 , J]} {\cal M}_{0,0}^{d_1 , d_2} (X,L ; J)$.
  The differential of the forgetful map $f_\bullet$ sends our basis $({\cal B}_1 , \stackrel{.}{r}_1 + r_1^* , r_2^* + \stackrel{.}{r}_2 ,  \nu_1 + \nu_2 , {\cal B}_2)$ 
 onto the family $({\cal B}_1 , r_1^* , r_2^*  ,  \nu_1 + \nu_2 , {\cal B}_2)$ of $T_{[u_1, u_2 , J]} {\cal M}_{0,0}^{d_1 , d_2} (X,L ; J)$.
 The canonical coorientation given by Lemma \ref{lemmawall} is defined such that $ \nu_1 $ is inward normal. The incidence index
 $\langle {\cal M}_{0,0}^{d_1, d_2} (X,L ; J) ,  {\cal M}_{0,0}^{(d_1, d_2)} 
(X,L ; J) \rangle$ thus equals $+1$ if and only if the basis  $( - \nu_1 , {\cal B}_1 , r_1^* , r_2^*  ,  \nu_1 + \nu_2 , {\cal B}_2)$
of $T_{[u_1, u_2 , J]} {\cal M}_{0,0}^{d_1 , d_2} (X,L ; J)$ is direct. The result thus follows from the parity of the cardinality of ${\cal B}_1$. $\square$

\section{Linking numbers, complete graphs and spanning subtrees}
\subsection{Complete graphs}
\subsubsection{Definitions}
\label{subsecdefs}

\begin{figure}[h]
\begin{center}
\includegraphics{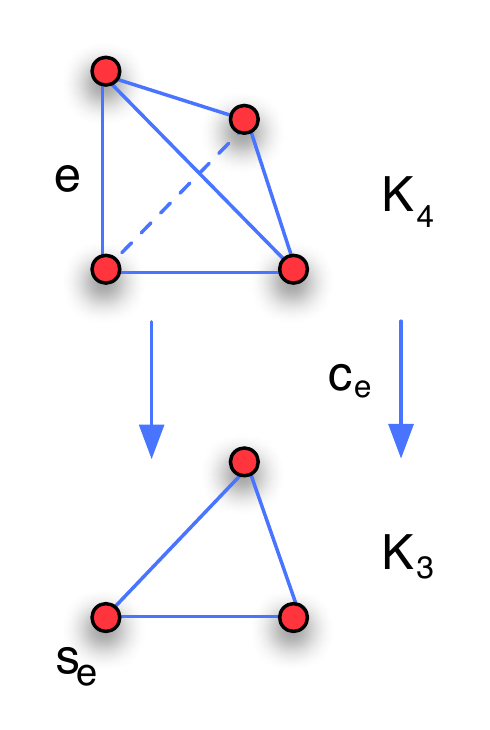}
\end{center}
\caption{Contraction map on $K_4$}
\label{figce}
\end{figure}

For every positive integer $n$, denote by $K_n$ the complete graph having $n$ vertices. Denote by $S_n$ (resp $E_n$) its set
of vertices (resp $n \choose {2}$ edges). For every $e \in E_{n+1}$, let $c_e : K_{n+1} \to K_n$ be the
contraction map of the edge $e$. It is defined as the quotient map of $K_{n+1} $ by the following equivalence relation ${\cal R}$ (see Figure
\ref{figce}):
\begin{itemize}
\item $\forall s_1 , s_2 \in S_{n+1} $, $s_1 {\cal R} s_2$ if and only if $s_1$ and $s_2$ bound $e$.
\item $\forall e_1 , e_2 \in E_{n+1} $, $e_1 {\cal R} e_2$ if and only if $e_1$ and $ e_2$ have one common boundary vertex and the other one
bounding $e$. 
\end{itemize}

The contraction map $c_e$ is surjective. The graph $K_n$ contains a unique vertex $s_e$ whose inverse image under $c_e$ does not reduce to a singleton,
but to the two bounding vertices of $e$. The inverse image to an edge of $K_n$ is made of one or two edges depending on whether this
edge contains $s_e$ in its boundary or not. 

Finally, for every commutative ring $A$ and every positive integer $n$, we denote by $A_n = \sum_{e \in E_n} A . e$ the free $A$-module
generated by the edges of $K_n$. For every edge $e_0$ of $K_{n+1}$, we then denote by $1_{e_0} = 1 . e_0$ the associated generator of the factor
$A . e_0$. The contraction map $c_{e_0} $ induces a morphism of $A$-modules 
$(c_{e_0} )_* : \sum_{e \inÊE_{n+1}} a_e . e \in A_{n+1} \mapsto \sum_{e \inÊE_{n+1}} a_e . c_{e_0} (e) \in A_{n} $ which contains
$A . e_0$ in its kernel.

\subsubsection{Spanning subtrees}

For every positive integer $n$, a subtree of $K_n$ is called a spanning subtree if and only if it contains all the vertices of $K_n$. 
These subtrees are maximal with respect to the inclusion, they form a set $T_n$ of cardinality $n^{n-2}$ as was established by
J. J. Sylvester and A. Cayley, see \cite{Syl}, \cite{Cay}. Likewise, for every $e \in E_n$, we denote by $T_n^e$ the subset of $T_n$ 
made of trees containing the edge $e$. 

\begin{lemma}
\label{lemmaTn}
Let $n \in \N^*$ and $e \in E_{n+1}$. Then, the contraction map $c_e$ induces a surjective map $T \in T_{n+1}^e \mapsto c_e (T) \in T_n$.
Moreover, the inverse image of an element $T \in T_n$ under this map contains $2^{\nu_T (s_e)}$ elements, where $\nu_T (s_e)$
denotes the valence of the special vertex $s_e$ in $T$. 
\end{lemma}

{\bf Proof:}

Let $T \in T_{n+1}^e$. Since $T$ spans $K_{n+1}$, $c_e (T) $ spans $ T_n$. Moreover, $c_e (T) $ cannot contain a cycle disjoint from $s_e$,
since the inverse image of such a cycle would produce a cycle in $T$. Likewise, it cannot contain a cycle containing $s_e$, since any
 inverse image of such a cycle in $K_{n+1}$ to which we add $e$ is a cycle of $K_{n+1}$. As a consequence, the image of every element
 of $T_{n+1}^e$ under $c_e$ is indeed a spanning subtree of $K_n$. 
 If now $T$ denotes an element of $T_n$, any inverse image of $T$ under $c_e$, to which we add $e$, provides an element of $T_{n+1}^e$,
 so that our map is surjective. The number of inverse images of $T$ is $2^{\nu_T (s_e)}$ since every edge of $K_n$ has one or two inverse images
 under $c_e$ depending on whether it contains $s_e$ or not in its boundary. $\square$\\

For every $T \in T_n$, denote by $E_n (T)$ the set of its $n-1$ edges and by $T_* :  \sum_{e \in E_n} a_e . e \in A_n \mapsto \prod_{e \in E_n (T)}
a_e \in A$ the associated $(n-1)$-linear form on the free $A$-module $A_n$ associated to any commutative ring $A$, see \S \ref{subsecdefs}.

\begin{defi}
\label{defforested}
For every positive integer $n$ and  commutative ring $A$, the $(n-1)$-linear form $\Phi_n = \sum_{T \in T_n} T_* : A_n \to A$ is called
the $n^{th}$ forested form. 
\end{defi}
Note that $K_1$ does not contain any edge, so that $A_1 = \{ 1 \}$ and that $\Phi_1 : A_1 \to A$ is the trivial map whose image is the unit
element of $A$. We now reach the aim of this paragraph, namely the following key Lemma \ref{lemmaforested}.

\begin{lemma}
\label{lemmaforested}
Let $n$ be a positive integer and $A$ a  commutative ring. Then, for every $a \in A_{n+1}$ and $e_0 \in E_{n+1}$,
$\Phi_{n+1} (a + 1_{e_0}) - \Phi_{n+1} (a) = \Phi_n \circ (c_{e_0} )_* (a)$, where $ \Phi_n$, $\Phi_{n+1}$ are the forested forms given
by Definition \ref{defforested}.
\end{lemma}

{\bf Proof:}

Let $T \in T_n$ and $T_1, \dots , T_{2^{\nu_T (s_e)}}$ be the inverse images of $T$ in $T_{n+1}^e$ given by Lemma \ref{lemmaTn}.
Then, for every $a \in A_{n+1}$, the composition $ T_* \circ (c_{e_0} )_* (a)$ writes
$$\prod_{e \in E_n (T) \, \vert \, s_e \notin \partial e} e_* \circ  (c_{e_0} )_* (a) \prod_{e \in E_n (T) \, \vert \, s_e \in \partial e} ((e_1)_* + (e_2)_* ) (a),$$
where for every $e \in E_n$, $e_*$ denotes the projection $a \in A_n \mapsto a_e \in A$ and if $s_{e_0}$ bounds $e$, $e_1$ and $e_2$ denote
the two inverse images of $e$ under $c_{e_0}$. By developing this expression, we deduce that 
$T_* \circ (c_{e_0} )_* (a) = \sum_{j=1}^{2^{\nu_T (s_e)}} \big( (T_j)_* (a + 1_{e_0}) - (T_j)_* (a) \big) \in A$. We then deduce from Lemma  \ref{lemmaTn}
after summation that $ \Phi_n \circ (c_{e_0} )_* (a) = \sum_{T \in T_{n+1}^{e_0}} (T_*  (a + 1_{e_0}) -T_*  (a ))  \in A$. The result now follows
from the fact that when $T \in T_{n+1} \setminus T_{n+1}^{e_0}$, the difference $T_*  (a + 1_{e_0}) -T_*  (a )$ vanishes by definition of $T_*$.
$\square$\\

Note that for every $n \in \N^*$,  the group of permutation of $S_n$ acts by permutation on $E_n$ and $T_n$.
It also acts by automorphisms on $A_n$, preserving the forested form $\Phi_n$.

\subsection{Linking numbers}
\label{subseclink}

Let $L$ be a closed oriented three-dimensional manifold and $A$ be a commutative ring. We denote by $C^0 (S^1 , L)$ the
space of continuous functions from the circle into $L$ and by
$C^0_A (S^1 , L) = \{ \gamma \in C^0 (S^1 , L) \; \vert \;  \gamma_* [S^1 ] = 0 \in H_1 (L ; A) \}$ those whose image vanish in homology. 
For every positive integer  $n$, we call link with $n$ components any element $( \gamma_1 , \dots ,  \gamma_n) \in C^0_A (S^1 , L)^n$ such that
the images of $ \gamma_i$, $i \in \{ 1 , \dots , n \}$, are disjoint from each other. We denote by $E^n (L , A)$ the space of links of $L$ having $n$ components
(homologically trivial in $H_1 (L ; A)$).

For every $( \gamma_1 ,  \gamma_2) \in E^2 (L , A)$, let $\Gamma_2$ be a two-chain of $L$ with coefficient in $A$ such that $\partial \Gamma_2 = \gamma_2$
as a one-cycle of $L$. We call linking number of $ \gamma_1 $ and $  \gamma_2$ the intersection index $ \gamma_1 \circ  \Gamma_2 \in A$ ; it is denoted
by $lk_2 ( \gamma_1 ,  \gamma_2)$. This index does not depend on the choice of $\Gamma_2$. It moreover satisfy the following properties.

{\bf $A_1$: Symmetry}

For every $( \gamma_1 ,  \gamma_2) \in E^2 (L , A)$, $lk_2 ( \gamma_1 ,  \gamma_2) = lk_2 ( \gamma_2 ,  \gamma_1)$.\\

{\bf $A_2$: Homotopy}

For every continuous function $f : [0,1] \to E^2 (L , A)$, $lk_2 \circ f :  [0,1] \to A$ is constant. \\

{\bf $A_3$: Vanishing}

For every contractile open subset $U$ of $L$ and every $( \gamma_1 ,  \gamma_2) \in E^2 (L , A)$ such that $\text{Im} (\gamma_1) \subset U$
and $\text{Im} (\gamma_2) \cap U = \emptyset$, $lk_2 ( \gamma_1 ,  \gamma_2) = 0$.\\

{\bf $A_4$: Additivity}

For every $( \gamma'_1 ,  \gamma_2) \in E^2 (L , A)$ and $( \gamma''_1 ,  \gamma_2) \in E^2 (L , A)$ such that $\gamma'_1 (1) = \gamma''_1 (1)$,
$lk_2 ( \gamma'_1 *   \gamma''_1 ,  \gamma_2) =  lk_2 ( \gamma'_1,  \gamma_2) + lk_2 (   \gamma''_1 ,  \gamma_2) $, where
$\gamma'_1 *   \gamma''_1$ denotes the concatenation of the  loops $ \gamma'_1 $ and $   \gamma''_1$.\\

{\bf $A_5$: Normalization}

Let $( \gamma_1 ,  \gamma_2) \in E^2 (L , A)$ be such that $ \gamma_2$ extends to en embedding of the unit disc $\Delta$ in $L$. 
Then, if the intersection index $ \gamma_1 \circ  \gamma_2 (\Delta)$ equals one, so does the linking number $lk_2 ( \gamma_1 ,  \gamma_2)$.\\

More generally, let $\underline{\gamma} = ( \gamma_1 , \dots ,  \gamma_n) \in E^n (L , A)$ be a link with $n$ components, $n \geq 1$. Let us choose
a bijection between the set $S_n$ of vertices of the complete graph $K_n$ and the set $\{ \gamma_1 , \dots ,  \gamma_n \}$. For every edge $e \in E_n$,
its boundary $\partial e$ consists of two vertices of $K_n$ which are then associated  to a link in $E^2 (L , A)$. We denote by $lk_2 (\partial e) \in A$ its
linking number. This provides an element $a_{\underline{\gamma}} = \sum_{e \in E_n} lk_2 (\partial e) . e $ in the $A$-module $ A_n$ introduced in
 \S \ref{subsecdefs} and we set $lk_n (\underline{\gamma} ) = \Phi_n (a_{\underline{\gamma}} ) \in A$, where $ \Phi_n$ is the forested form given
 by Definition \ref{defforested}. This element $lk_n (\underline{\gamma} ) $ does not depend on the choice of the bijection 
 $S_n \to \{ \gamma_1 , \dots ,  \gamma_n \}$.

\begin{defi}
\label{defselflinkingweight}
For every link $\underline{\gamma} \in E^n (L , A)$ with $n$ components, $n \geq 1$, the element $lk_n (\underline{\gamma} )  \in A$ is called its
self-linking weight. 
\end{defi}
In particular, when $n=1$, $lk_1 (\gamma) = 1$ whatever $\gamma \in E^1 (L , A)$ is. For $n = 2$, the self linking weight  given by Definition
 \ref{defselflinkingweight} is just the linking number of the two components of the link. The self-linking weight satisfies the following properties. 
 
 {\bf $B_1$: Symmetry}

For every $( \gamma_1 , \dots ,  \gamma_n) \in E^n (L , A)$ and every permutation $\sigma $ of $\{ 1 , \dots , n \}$, $n \geq 1$, 
$lk_n ( \gamma_1 , \dots ,  \gamma_n) = lk_n ( \gamma_{\sigma (1)} , \dots ,  \gamma_{\sigma (n)})$.\\

{\bf $B_2$: Homotopy}

For every continuous function $f : [0,1] \to E^n (L , A)$, $lk_n \circ f :  [0,1] \to A$ is constant. \\

{\bf $B_3$: Vanishing}

For every contractile open subset $U$ of $L$ and every $( \gamma_1 , \dots , \gamma_n) \in E^n (L , A)$ such that $\text{Im} (\gamma_1) \subset U$
and $\text{Im} (\gamma_j) \cap U = \emptyset$, $j \in \{ 2 , \dots , n \}$,  $lk_n ( \gamma_1 ,  \dots , \gamma_n) = 0$.

\section{Open Gromov-Witten invariants in dimension six}

\subsection{Statement of the result}

Let $(X , \omega)$ be a closed symplectic manifold of dimension six.
Let $L \subset X$ be a closed orientable Lagrangian submanifold equipped with a Spin structure (and thus an orientation).
We assume that there is a commutative ring $A$ such that the inclusion $L \to X$ induces an injective morphism
$H_1 (L ; A) \to H_1 (X ; A)$ in homology. We denote by ${\cal J}_\omega$ the space of almost-complex structures of $X$
tamed by $\omega$ and of class $C^l$, $l \gg 1$. We then denote by ${\cal J}_\omega^{\neq 0}$ a connected open subset of
${\cal J}_\omega$ such that for every $J \in {\cal J}_\omega^{\neq 0}$,  $X$ contains no Maslov zero
$J$-holomorphic disc with boundary on $L$. 

Let $d \in H_2 (X , L ; \Z)$ be of positive Maslov index $\mu_L (d)$, $r,s$ be non-negative integers and $J \in {\cal J}_\omega^{\neq 0}$
be generic. The space $ {\cal M}_{r,s}^d (X,L ; J) $ of simple $J$-holomorphic discs with boundary on $L$, homologous to $d$
and having $r$ (resp. $s$) marked points on their boundaries (resp. interiors) is an oriented manifold of dimension $\mu_L (d) + r +2s$
whose compactification has boundaries and corners, see \S \ref{subsecsimple}. It is equipped with an evaluation map at the marked points
denoted by $ev : {\cal M}_{r,s}^d (X,L ; J) \to L^r \times X^s$. Likewise, for every positive integer $n$ and every $d_1 , \dots , d_n \in H_2 (X, L ; \Z)$ such that 
$\mu_L (d_i) >0$ and $d_1 + \dots + d_n = d$, the moduli space ${\cal M}_{r,s}^{d_1, \dots , d_n} (X,L ; J) $ of simple reducible 
$J$-holomorphic discs with $n$ components homologous to $d_1 , \dots , d_n$ and having $r$ (resp. $s$) marked points on their boundaries 
(resp. interiors) is an oriented manifold of dimension $\mu_L (d) + r +2s$ whose compactification has boundaries and corners, 
see \S \ref{subsecreducible}.  It is equipped with an evaluation map at the marked points denoted by
$ev : {\cal M}_{r,s}^{d_1, \dots , d_n} (X,L ; J) \to L^r \times X^s$.

  Let us denote by $\stackrel{\! \!\!\!\!\!\! \!\!\!\!\!\!\! \! \!\!\!\!\! \circ}{{\cal M}^{d_1, \dots , d_n}_{r,s}} (X,L ; J) $ the dense open subset of ${\cal M}_{r,s}^{d_1, \dots , d_n} (X,L ; J) $ 
  made of discs whose $n$ components have pairwise disjoint boundaries in $L$. It is equipped with a boundary map
  $\partial : D \in   \stackrel{\! \!\!\!\!\!\! \!\!\!\!\!\!\! \! \!\!\!\!\! \circ}{{\cal M}^{d_1, \dots , d_n}_{r,s}} (X,L ; J) \mapsto \partial D \in E^n (L , A)$ in the space of $n$-components links
  of $L$, see \S \ref{subseclink}. Every component of such a link indeed bounds in $X$ and thus under our hypothesis also bounds 
  a two-chain in $L$. This open subset thus gets equipped with a locally constant function of self-linking weight
  $lk_n : D \in \stackrel{\! \!\!\!\!\!\! \!\!\!\!\!\!\! \! \!\!\!\!\! \circ}{{\cal M}^{d_1, \dots , d_n}_{r,s}} (X,L ; J) \mapsto lk_n (\partial D) \in A$, see Definition \ref{defselflinkingweight}.

Let finally 
$[ {\cal M}_{d, r,s} (X,L ; J) ] = \sum_{n=1}^\infty \frac{1}{n!} \sum_{d_1 + \dots + d_n = d} lk_n [\stackrel{\! \!\!\!\!\!\! \!\!\!\!\!\!\! \! \!\!\!\!\! \circ}{{\cal M}^{d_1, \dots , d_n}_{r,s}} (X,L ; J)]$
be the fundamental class of the $CW$-complex $ {\cal M}_{d, r,s} (X,L ; J) $ and
$$ev_* [ {\cal M}_{d, r,s} (X,L ; J) ] = \sum_{n=1}^\infty \frac{1}{n!} \sum_{d_1 + \dots + d_n = d} lk_n 
ev_*[\stackrel{\! \!\!\!\!\!\! \!\!\!\!\!\!\! \! \!\!\!\!\! \circ}{{\cal M}^{d_1, \dots , d_n}_{r,s}} (X,L ; J)]$$
be the image of this fundamental class under the evaluation map in the space of chains of $L^r \times X^s$ with coefficients in $A$ and dimension
$\mu_L (d) + r +2s$. Note that the sums in the right hand sides of these equalities are actually finite and that the term $\frac{1}{n!}$
compensates the order of the $n$-tuples $(d_1, \dots , d_n)$, since the group of permutation of this set acts by preserving the fundamental classes. 

\begin{theo}
\label{theoGW}
Let $(X , \omega)$ be a closed symplectic manifold of dimension six.
Let $L \subset X$ be a closed orientable Lagrangian submanifold equipped with a Spin structure.
Let $A$ be a commutative ring such that the inclusion $L \to X$ induces an injective morphism
$H_1 (L ; A) \to H_1 (X ; A)$. Let ${\cal J}_\omega^{\neq 0}$ be a connected open subset of
the space of almost-complex structures tamed by $ \omega$ such that for every $J \in {\cal J}_\omega^{\neq 0}$,  $X$ contains no Maslov zero
$J$-holomorphic disc with boundary on $L$. Then, for every $d \in H_2 (X , L ; \Z)$ of positive Maslov index,
every $r \geq 1$, $s \geq 0$ and $J \in {\cal J}_\omega^{\neq 0}$ generic, the chain
$ev_* [ {\cal M}_{d, r,s} (X,L ; J) ]$ is a cycle whose homology class in $H_{\mu_L (d) + r +2s} (L^r \times X^s ; A)$ does not depend on
the generic choice of $J$.
\end{theo}

When $r=0$, Theorem \ref{theoGW} does not hold in general, since the  $CW$-complex $ {\cal M}_{d, r,s} (X,L ; J) $ might then have an extra
boundary component not contracted under the evaluation map. This extra component parameterizes discs whose boundary get shrunk to a point,
that is of $J$-holomorphic spheres meeting $L$. We restrict ourselves to the open subset $ {\cal J}_\omega^{\neq 0}$
in order not to have to take into account branched covers of Maslov zero pseudo-holomorphic discs. 
Note that when $L$ is a Lagrangian sphere in a Calabi-Yau manifold $X$, there exists for every energy bound $E$ a generic almost-complex
structure $J$ tamed by $\omega$ such that $X$ does not contain any  $J$-holomorphic disc with boundary on $L$ and energy less than $E$,
see Theorem $1.6$ of \cite{WelsOpt} (or also Corollary $4.3$ of \cite{WelsFloer} , Theorem $4.1$ of \cite{WelsICM}).
Finally, the hypothesis on the ring $A$ makes it possible to define the self-linking weights $(lk_n)_{n \geq 1}$
and thus to kill the boundary components of ${\cal M}_{r,s}^d (X,L ; J)$.

\begin{cor}
\label{corGW}
Under the hypothesis of Theorem \ref{theoGW}, for every classes $\alpha_1, \dots , \alpha_r \in H^* (L ; A)$ and 
$\beta_1, \dots , \beta_s \in H^* (X ; A)$ such that $\sum_{i=1}^r \deg (\alpha_i)+ \sum_{j=1}^s \deg (\beta_j) = \mu_L (d) + r +2s$
and $\sum_{i=1}^r \deg (\alpha_i) >0$, the open Gromov-Witten invariant 
$$GW_d (X,L ; \alpha_1, \dots , \alpha_r ,  \beta_1, \dots , \beta_s) =
 \int_{ev_* [ {\cal M}_{d, r,s} (X,L ; J) ] } \alpha_1\wedge \dots \wedge \alpha_r
\wedge \beta_1\wedge \dots \wedge \beta_s \in A$$
does not depend on the generic choice of $J  \in {\cal J}_\omega^{\neq 0}$. $\square$
\end{cor}

Given generic submanifolds $A_1 , \dots , A_r$ (resp. $B_1, \dots , B_s$) of $L$ (resp. $X$) Poincar\'e duals to $\alpha_1, \dots , \alpha_r $
(resp. $\beta_1, \dots , \beta_s $), the invariant $GW_d (X,L ; \alpha_1, \dots , \alpha_r ,  \beta_1, \dots , \beta_s)$ given by Corollary \ref{corGW} counts the number of
$J$-holomorphic discs which are either irreducible or reducible not connected, have boundary on $L$ with total relative homology class $d$
and which meet $A_1 , \dots , A_r$ as well as $B_1, \dots , B_s$. These discs are counted with respect to some sign and are weighted by
the self-linking weights of their boundaries. 

\subsection{Proof of Theorem \ref{theoGW}}

Let $J  \in {\cal J}_\omega^{\neq 0}$ be a generic almost-complex structure. The boundary of the chain $ev_* [ {\cal M}_{d, r,s} (X,L ; J) ]$ is
carried by the image under the evaluation maps of codimension one components  of the boundary $\partial {\cal M}_{d, r,s} (X,L ; J)$
of the moduli space space ${\cal M}_{d, r,s} (X,L ; J)$ in its stable maps compactification, see \cite{FOOO}. From Gromov's compactness theorem
(see \cite{Frauen}), these codimension one components are themselves moduli spaces of nodal $J$-holomorphic discs with boundary on $L$.
These nodal discs may contain a priori any number of nodes as well as spherical components attached to them. Moreover, some of these
spherical component might be multiply covered while the discs might not be simple in the sense of Definition \ref{defsimple}.

We prove in \S \ref{subsubsectheosimple} that only simple $J$-holomorphic discs may appear in codimension one boundary components
of $ev_* [ {\cal M}_{d, r,s} (X,L ; J) ]$, thanks to the hypothesis made on the absence of Maslov zero
disc. We then prove in \SÊ\ref{subsubsecspheres} thanks to the hypothesis $r \geq 1$ that no
spherical component may appear in codimension one boundary components
of $ev_* [ {\cal M}_{d, r,s} (X,L ; J) ]$, so that the only codimension one boundary components that may appear parameterize nodal reducible
 $J$-holomorphic discs having a unique node. Finally, we prove in \S \ref{subsubseccontrib} that each of the latter codimension one boundary components 
 have a vanishing weight in the boundary of the chain $ev_* [ {\cal M}_{d, r,s} (X,L ; J) ]$, which follows from the key property of the self-linking weight
 given by Lemma \ref{lemmaforested}. It follows that $ev_* [ {\cal M}_{d, r,s} (X,L ; J) ]$ is a cycle. Let now $J_0 , J_1$ be generic elements
 of ${\cal J}_\omega^{\neq 0}$ and $(J_t)_{t \in [0,1]}$ be a generic path joining $J_0$ to $J_1$. It follows exactly along the same lines that
 the boundary of the chain $ev_* [ {\cal M}_{d, r,s} (X,L ; \cup_{t \in [0,1]} J_t) ]$ reduces to the difference 
 $ev_* [ {\cal M}_{d, r,s} (X,L ; J_1) ] - ev_* [ {\cal M}_{d, r,s} (X,L ; J_0) ]$ in the space of chains of dimension $\mu_L (d) + r +2s$ in $L^r \times X^s$.
 Theorem \ref{theoGW} follows.
 
 \subsubsection{Theorem of decomposition into simple discs}
\label{subsubsectheosimple} 

We recall in this paragraph the theorem of decomposition into simple discs established by Kwon-Oh and Lazzarini, see  \cite{KwonOh}, \cite{Laz}, \cite{Laz1}. 
\begin{theo}
\label{theodecomp}
Let $L$ be a closed Lagrangian submanifold of a $2n$-dimensional closed symplectic manifold $(X , \omega)$. Let $u : (\Delta , \partial \Delta) \to (X,L)$
be a non-constant pseudo-holomorphic map. Then, there exists a graph ${\cal G} (u)$ embedded in $\Delta$ such that $\Delta \setminus {\cal G} (u)$
has only finitely many connected components. Moreover, for every connected component ${\cal D} \subset \Delta \setminus {\cal G} (u)$, there exists a
surjective map $\pi_{\overline{\cal D}} : \overline{\cal D} \to \Delta$, holomorphic on ${\cal D}$ and continuous on $\overline{\cal D}$, as well as a simple
pseudo-holomorphic map $u_{\cal D} : \Delta \to X$ such that $u\vert_{\overline{\cal D}} = u_{\cal D} \circ \pi_{\overline{\cal D}}$. The map $\pi_{\overline{\cal D}}$
has a well defined degree $m_{\cal D} \in \N$, so that $u_* [\Delta] = \sum_{\cal D} m_{\cal D} (u_{\cal D})_* [\Delta] \in H_2 (X,L ; \Z)$, the sum being taken
over all connected components ${\cal D}$ of $ \Delta \setminus {\cal G} (u)$. $\square$
\end{theo}
The graph $ {\cal G} (u)$ given by Theorem \ref{theodecomp} is called the frame or non-injectivity graph, see \cite{Laz}, \cite{Laz1} (or \S $3.2$ of \cite{BirCor}) for its definition. 

Here, the Lagrangian $L$ being supposed to be orientable, all Maslov indices of discs with boundaries on $L$ are even. Moreover, since $J$ is generic,
none of these indices are negative by Theorem \ref{theodim} and then since $J  \in {\cal J}_\omega^{\neq 0}$, they are all greater or equal to two. As a consequence,
as soon as a degree $m_{\cal D}$ given by Theorem \ref{theodecomp} is greater than one, the codimension of the image of such discs under the
evaluation map in the chain $ev_* [ {\cal M}_{d, r,s} (X,L ; J) ]$ is at least two, since this image is carried by simple discs of total Maslov index
bounded from above by $\mu_L (d) -2$. Moreover, the non-injectivity graphs given by Theorem \ref{theodecomp} are then trivial in codimension one,
since a pair of simple discs only meet at one point in codimension one, $J$ being generic. Hence, the only discs to consider in order to prove
Theorem \ref{theoGW} are simple discs. 

 \subsubsection{Pseudo-holomorphic spheres}
 \label{subsubsecspheres} 
 
 The condition for a simple $J$-holomorphic sphere to meet a non-constant $J$-holomorphic disc costs two degrees of freedom.
 Moreover, the Maslov index of such a $J$-holomorphic sphere can only increase under branched coverings, $J$ being generic, so that
 the same holds for their Fredholm index. As a consequence,  the moduli spaces of stable maps containing spherical components are
 all of codimension at least two in the compactification of $ {\cal M}_{d, r,s} (X,L ; J) $ and thus do not contribute to the boundary 
 $ev_* [ {\cal M}_{d, r,s} (X,L ; J) ]$.
 
 Now, the condition for a simple $J$-holomorphic sphere to meet a constant $J$-holomorphic disc costs only one degree of freedom, since such
 a disc lies in $L$ and the moduli space of $J$-holomorphic spheres meeting $L$ is of codimension one in the moduli space of homologous 
$J$-holomorphic spheres. Likewise, the moduli space of reducible $J$-holomorphic discs with boundary on $L$, homologous to $d$ and 
having exactly one spherical component attached to a ghost component is of dimension $\mu_L (d) -1$ when it has no marked point. If the disc is irreducible, 
all the points marked on the boundary of the disc evaluate to the same image point, so that this stratum is of codimension $r+1$ in the compactification 
of $ {\cal M}_{d, r,s} (X,L ; J) $. Since we assumed that $r \geq 1$, this stratum does not contribute to the boundary of $ev_* [ {\cal M}_{d, r,s} (X,L ; J) ]$.

If the disc is reducible, the boundary of the ghost component might contain none of the $r$ marked points. 
However, from the vanishing property $B3$ of the self-linking weight given in \S \ref{subseclink}, all self-linking weights of these reducible discs
appearing in codimension one vanish. As a consequence, their contribution to the boundary of the chain $ev_* [ {\cal M}_{d, r,s} (X,L ; J) ]$ also vanish.

\subsubsection{Contribution of nodal discs}
\label{subsubseccontrib}

In order to prove Theorem  \ref{theoGW}, we thus now just have to consider the contribution of nodal discs with a unique node to
the  boundary of the chain $ev_* [ {\cal M}_{d, r,s} (X,L ; J) ]$. These are the discs studied in \S \S \ref{subsecred} and \ref{subsecreducible}.
Let $n \geq 1$ and ${\cal M}_{r,s}^{(d_1, d_2), d_3, \dots , d_{n+1}} (X,L ; J) $ be such a moduli space equipped with its orientation. 
From Lemma \ref{lemmawall} we know that this space, or rather its image under the forgetful map $f_\bullet$, is of codimension one
in ${\cal M}_{r,s}^{d_1, d_2, d_3, \dots , d_{n+1}} (X,L ; J) $ and canonically cooriented. Let  
$[u_1^\lambda , u_2^\lambda, u_3 , \dots , u_{n+1}, J , \underline{z} , \underline{\zeta}]_{\lambda \in ]-\epsilon, \epsilon[}$ be a path 
of $ {\cal M}_{r,s}^{d_1, \dots , d_{n+1}} (X,L ; J)$ positively transversal to the wall ${\cal M}_{r,s}^{(d_1, d_2), d_3, \dots , d_{n+1}} (X,L ; J) $ at $\lambda = 0$.
By definition of the canonical coorientation given by Lemma \ref{lemmawall} and from the property of normalization $A5$ 
given in \S \ref{subseclink}, the linking number $lk_2 (u_1^\lambda (\partial \Delta) , u_2^\lambda (\partial \Delta))$ between the boundaries of the
two first discs increase by one while crossing the parameter $\lambda = 0$. 

Let us denote by $a^- \in A_{n+1} $ (resp. $a^+ \in A_{n+1} $) the element
of $A_{n+1} $ associated to the link $(u_1^\lambda (\partial \Delta), u_2^\lambda (\partial \Delta), u_3 (\partial \Delta), \dots , u_{n+1} (\partial \Delta))$
for $\lambda < 0$ (resp. $\lambda > 0$), so that its self-linking weight $lk_{n+1} (u_1^\lambda (\partial \Delta), u_2^\lambda (\partial \Delta), u_3 (\partial \Delta), \dots , u_{n+1} (\partial \Delta))$ 
equals by definition $\Phi_{n+1} (a^-)$ (resp. $\Phi_{n+1} (a^+)$), see \S \ref{subseclink}. By definition,
$a^+ = a^- + 1_e$, where $e$ denotes the edge of $K_{n+1}$ joining the vertices associated to $u_1^\lambda (\partial \Delta)$ and 
$ u_2^\lambda (\partial \Delta)$. From Proposition \ref{propboundaries} follows that the chain $[{\cal M}_{r,s}^{(d_1, d_2), d_3, \dots , d_{n+1}} (X,L ; J) ]$
contributes to the boundary of $lk_{n+1} [\stackrel{\! \!\!\!\!\!\!\!\!\!\!\!\!\! \!\!\!\!\!\!\!\!\!\!\!\!\!  \circ}{{\cal M}^{d_1, \dots , d_{n+1}}_{r,s}}  (X,L ; J)]$
 with the coefficient $\Phi_{n+1} (a^-  + 1_e) - \Phi_{n+1} (a^-)$. 

Now, from the same Proposition \ref{propboundaries} and the additivity property $A4$ given in \S \ref{subseclink}, 
the same chain $[{\cal M}_{r,s}^{(d_1, d_2), d_3, \dots , d_{n+1}} (X,L ; J) ]$
contributes to the boundary of 
$lk_{n} [\stackrel{\! \!\!\!\!\!\!\!\!\!\!\!\! \!\!\!\!\!\!\! \!\!\!\!\!\!\!  \!\!\!\!\!\!\! \!\!\!\!\!\!\! \!\!\!\! \circ}{{\cal M}^{d_1+d_2 , d_3, \dots , d_{n+1}}_{r,s}} (X,L ; J)]$
with the coefficient
$-\Phi_n \circ (c_e)_* (a^-)$. From Lemma \ref{lemmaforested} we conclude by summation that the contribution of this chain 
to the boundary of $ [ {\cal M}_{d, r,s} (X,L ; J) ]$ vanishes. Hence the result. $\square$

\addcontentsline{toc}{part}{\hspace*{\indentation}Bibliography}

\bibliographystyle{abbrv}

\vspace{0.7cm}
\noindent 
Universit\'e de Lyon ; CNRS ;
Universit\'e Lyon~1 ; Institut Camille Jordan

\end{document}